\documentclass[12pt]{article}
\usepackage[T1]{fontenc}
\usepackage{multicol,epsfig,amssymb,amsmath,mathrsfs,verbatim} 
\newcommand{\proof}{\medskip \noindent {\bf Proof. }}
\newcommand{\qed}{\null\hfill  $\Box\;\;$ \medskip}
\newcommand{\QED}{\null\hfill \hskip 1cm\Box\;\; \medskip}

\usepackage{xcolor}

\newtheorem {theorem} {Theorem}

\newtheorem {definition} {Definition}
\newtheorem {lemma} {Lemma}
\newtheorem {remark} {Remark}

\begin{document}
\def\R{\mathbb{R}}
\def\C{\mathbb{C}}
\def\N{\mathbb{N}}
\def\Z{\mathbb{Z}}
\def\Q{\mathbb{Q}}
\def\Ra{\operatorname{\textrm{rank}}}
\def\sign{\operatorname{sign}}
\def\ker{\operatorname{ker}}
\def\x{\times}
\def\<{\langle}
\def\>{\rangle}
\title{On a 4-dimensional  subalgebra of the 12-tone Equal Tempered Tuning}
\author{J\'an Halu\v{s}ka$^1$, Ma\l{}gorzata Jastrz\k{e}bska$^2$}
\date{}
\maketitle
\begin{abstract}
An operation of  associative, commutative and distributive multiplication on {   Euclidean vector space} $\mathbb{E}_4$  is introduced by a skew circulant  matrix.  The resulting algebra $\mathbb{W}$ over $\mathbb{R}$  is isomorphic to $\mathbb{C} \times \mathbb{C}.$  The related algebraic, geometrical, and topological properties are  given.There are  subplanes of $\mathbb{W}$ isomorphic to the Gauss and Clifford complex number planes. A topology  on $\mathbb{W}$  is given by a norm which is a sum of two norms. A hint how to apply this 4 dimensional algebra over $\mathbb{R}$  to the 12-tone Equally Tempered Tuning algebra is given. 
\end{abstract}

	
	
\noindent {\it Acknowledgement.} The paper is supported with the Grant  VEGA~2/0106/19 (Wooden pipe configuration of  historic organ positives in Slovakia).

\noindent{\it Author's Address:} $^1$ Math. Inst., Slovak Acad. Sci., Gre\v{s}\'{a}k Str. 6, Ko\v{s}ice, Slovakia 
                                  
                                  $^2$  Siedlce University,  Faculty of Exact and Natural Sciences, 08-110 Siedlce, ul. 3 Maja 54, Poland

\noindent{\it e-mail:} $^1$ jhaluska @ saske.sk, 
 $^2$ majastrz2 @ wp.pl
\section{Musical motivation and introduction} 

The {\it 12-tone Equal Temperament Tuning}, 12-TET, of musical instruments is generally and widely used today in music. It is arose from the  European aesthetic models of harmony; from the  structure of  overtones which are automatically physically  generated in the process of vibration of every material string (or, equivalently, of vibrations of a  column of air in every  pipe); and of cause,  from the mathematical number theory discoveries.   A real multiplicative number interval $[1, 2)$ is multiplicatively divided by number 12, so there are 12 tones with frequencies $$\{ a_0 (\root{12}\of {2})^{0 + 12 k}, a_0 (\root{12}\of{2})^{1+ 12 k}, \dots, a_0 (\root{12}\of{2})^{p+ 12 k} \},$$ $p = 0,1,2, \dots, 11$. This group of tones is  called an {\sl ($k-$) octave} (in the musical theory). There are used 8 octaves in music, the lower and upper octaves are not complete 12 in the audio range of human ear in the standard keyboard of clavier instruments. The octave names are: Subcontra, Contra, Great, Small, One-line, Two-line, Three-line, Four-line, (theoretically, Five-line, etc.). The value $a_0$ is usually a {\it camertone} which has the frequency of sound 440~Hz today.

A mathematics  in this  12-tone Equal Temperament Tuning is locally similar to a   structure of a vector algebra in the middle area of the audio diapason. The two vector operations have the following  musical sense.  {\it Melody} is bounded  with a role of vector addition, while {\it harmony}  is bounded with the operation  of multiplication of vectors (tones). An operation  expressing the {\it loudness} of the tone plays a role of the mixed operation of multiplication vectors by scalars. Thus, every European tonal music can be  mathematically understood as a special coding of pictures in an algebra (over $\mathbb{R}$) which appear in the space of 12-TET tone (vector) space. It is an algebra over the real line acting in time. In other words, it is a kind of  "movie" which we can hear (not to see). 

Musical sound (= tone) is mathematically modelled for a string (equivalently, a column of air in the pipe)  as a  Fourier series.  An operation of addition  of tones is  reflected  as a sum of tones. It is  good observable when combining organ registers.  

The definition of the operation of multiplication in music  is defined  by  a  special way, which seems rather artificial from a superficial view. Not going into details,
the definition of it is commonly  know among musicians  as a {\sl spiral of fifths}. Mathematically, an operation of multiplication of vectors  is defined by a skew circulated matrices.

Let us note  the following four psycho-acoustical phenomena which are present/asked from the multiplication operation:

1. The task of freely transposition, every melody may be be started from an arbitrary tone; Temperaments which are not 12-TET do not have this property. 

2. All tones which have a octave-shifted frequency are  mentioned harmonically to be equal; this musical phenomenon  is called  the {\sl octave equivalence}.  

3.  A human ear do not psycho-physically  distinguish individual tones $T$ and $(-T)$ although there exists a real physical situation when $T \oplus (-T)=0$. For instance, two  sounding  equal organ pipes with opposite standing cutout slits  produce together physically  a zero  sound (the two sound  waves annihilate one other). The loudness of tones (multiplication by scalar) is given by amplitude of tone vibrations, it is denoted as a usual multiplication by a real number. 
 
4. In the case of the 12-tone Equal Temperament, the  spiral of fifths is approximated (deformed) to the 
{\sl circle of fifth}, all semitones are equal, hence the name of the tuning.

 These phenomena are taken into account 
when introducing the definition of the operation of multiplication in  the algebra. To compare these requirements with  the resulting formalized  mathematical working definition, cf. Definition \ref{D1}.

For a large explanatory mathematical and also music-theoretical context, cf. a book \cite{Haluska}. From this book  it follows  that the all 12-tone Equal Temperament music is  about an humanly acceptable appropriate psycho-acoustical mathematical approximation of musical sound.

When considering every tonal music as a movie of  objects  we need to know the arithmetic in this 12-dimensional algebra over $\mathbb{R}$. Firstly, we are able to describe the  algebra of the dimension 4 which is a  subalgebra of the 12-TET algebra (in a musical terminology, the basic tones are $c, dis, fis, a$). The further subalgebra of the 12-TET are the whole-tone algebra (in a musical terminology, the basic tones are $c, d, e, fis, gis, ais$), to be complete, there is also a triton subalgebra  of 12-TET (in a musical terminology, the basic tones are$c, fis$) which has rather no importance in music.

In the accord to the previous motivations, an  4-dimensional algebra $\mathbb{W}$    over $\mathbb{R}$  (a subalgebra of 12-TET, isomorphically described) is investigated in this paper. We give the related algebraic, geometrical, and topological properties of this algebra.  In Chapter 2 we introduce associative, commutative and distributive multiplication on $\mathbb{E}_4$. As in a Theorem \ref{1MJ} this algebra can be seen as a generalization of complex numbers.

A well known extension of the set of Gaussian complex numbers to four dimensions is a algebra of quaternions described by W. R. Hamilton in 1843, \cite{Hamilton}.  A feature of quaternions in general is that they form a noncommutative division algebra over $\mathbb{R}$.  In \cite{Tosun},   the subset of commutative quaternions  was studied.  Unfortunately, we cannot use the subalgebra of commutative quaternions because a skew circular structure  is needed  by musicians.

All definitions and notions about algebras over fields used in this paper one can find for example in \cite{DK}.  We use Hankel, Toeplitz and skew circulant matrices (see e.g.~\cite{Davis}).   
Recall that  $n\times n$ matrix 
$T=(\tau_{j, k})^{n-1}_{j, k=0}$ is said to be {\it Toeplitz} if each descending diagonal from left to right is constant i.e.,
 $\tau_{j, k}= \tau_{j - k}$.   The 
 skew-circulant matrix is a special cases of Toeplitz matrix.
A  $n\times n$ matrix   $S=(\sigma_{j, k})^{n - 1}_{j, k=0}$ is said to be skew-circulant 
if $\sigma_{j, k}=\sigma_{j - k}$ and $\sigma_{- l}= - \sigma_{n - l}$ for $1 \leq l \leq n - 1$.
Skew circulant matrices have  applications to various disciplines (see ~\cite{Boettcher, LC, Ng,  YY, Yao}, and the references given there).
A square matrix in which each ascending skew-diagonal from left to right is constant is know as Hankel matrix.

In Chapter 2, in order to describe our 4-dimensional algebra  $\mathbb{W}$  over $\mathbb{R}$,  we give  its isomorphisms to known mathematical objects.
In Theorem \ref{rep1} we indicate matrix representation of elements $\mathbb{W}$ where we use  skew circulant matrices. 
In Theorem \ref{TH8} it is shown that the algebra $\mathbb{W}$  is isomorphic to $\mathbb{C} \times \mathbb{C}.$
In Theorem \ref{TH10} we indicate the method of finding invertible elements and its inverse.
Subplanes in $\mathbb{W}$ isomorphic to the Gaussian and Clifford complex numbers are indicated in  Chapter~\ref{Ch3} and Chapter~\ref{Ch4}.

\section{Operation of multiplication. \\ 
Isomorphisms of algebras.}
Let us consider the Euclidean vector space  $\mathbb{E}_4$ over $\mathbb{R}$  together with the standard vector addition, denoted $\oplus,$ and scalar multiplication. So that each element $\mathbf{x} \in \mathbb{E}_4$ 
  is expressible
  $$\mathbf{x} =  X_{\mathbf{1}} \  {\mathbf{1}}  
  \oplus  X_\mathbf{i}\ \mathbf{i}\oplus X_\mathbf{j} \ \mathbf{j} \oplus  X_\mathbf{k} \ \mathbf{k},$$  
   as a linear combination of  the standard basis elements,
$${\mathbf{1}} = (1,0,0,0), \mathbf{i}=(0,1, 0, 0), \mathbf{j}=(0, 0, 1, 0),  \mathbf{k}=(0, 0, 0,  1),$$ 
where $X_{\mathbf{1}}, X_\mathbf{i}, X_\mathbf{j}, X_\mathbf{k} \in   \mathbb{R}$. 

We use also the notation:
$\mathbf{x} =  (X_{\mathbf{1}},  X_\mathbf{i},  X_\mathbf{j}, X_\mathbf{k})$ and $\Lambda=(0,0,0,0)$. 
A~sign $\ominus$ denotes the inverse group operation to $\oplus$ in  $\mathbb{E}_4$.   

Let us define an operation of multiplication in $\mathbb{E}_4$ as follows: 
\begin{definition}\rm \label{D1} 
	 Let $\mathbf{x} = (X_\mathbf{1}, X_\mathbf{i}, X_\mathbf{j}, X_\mathbf{k})$, $\mathbf{y} = (Y_{\mathbf{1}
	 }, Y_\mathbf{i}, Y_\mathbf{j}, Y_\mathbf{k}) \in \mathbb{E}_4$. Then
	$$\mathbf{x} \otimes \mathbf{y} \stackrel{def}{=}
(X_{\mathbf{1}
} Y_{\mathbf{1}
}  - X_\mathbf{\mathbf{i}} Y_\mathbf{k}  -  X_\mathbf{j} Y_\mathbf{j} - X_\mathbf{k} Y_\mathbf{i}) {\mathbf{1}
}\oplus (X_{\mathbf{1}
} Y_\mathbf{i}+ X_\mathbf{i}Y_\mathbf{1}  -X_\mathbf{j} Y_\mathbf{k} -  X_\mathbf{k} Y_\mathbf{j} ) \mathbf{i}$$
$$\oplus ( X_{\mathbf{1}
} Y_\mathbf{j} + X_\mathbf{i}Y_\mathbf{i}+ X_\mathbf{j} Y_{\mathbf{1}
}    - X_\mathbf{k} Y_\mathbf{k})\mathbf{j}
\oplus (X_{\mathbf{1}
} Y_\mathbf{k}  +  X_\mathbf{i}Y_\mathbf{j}  + X_\mathbf{j} Y_\mathbf{i}+ X_\mathbf{k} Y_{\mathbf{1}
}) \mathbf{k}  \in \mathbb{E}_4.$$ 

Let us denote by ~$\mathbb{W}$ the space $\mathbb{E}_4$ equipped with the operation  of multiplication $\otimes$.

 \end{definition}	
	
The proof of the following theorem is trivial and therefore omitted. 

\begin{theorem}  The set  $\mathbb{W}$ is an associative, commutative algebra over a field $\mathbb{R}$  with the  multiplicative identity $\mathbf{1}=(1,0,0,0)$.
\end{theorem}
 
\begin{remark}\rm  
A topology on  $\mathbb{W}$ is given with a standard Euclidean  norm $$|||\mathbf{x}||| = \sqrt{X_\mathbf{1}^2 + X_\mathbf{i}^2 + X_\mathbf{j}^2 + X_\mathbf{k}^2} = \sqrt { <\mathbf{x}, \mathbf{x}>},$$  implied  from the scalar product
$<\mathbf{x},\mathbf{y}> = \frac{|||\mathbf{x} \oplus \mathbf{y}|||^2 - |||\mathbf{x} \ominus  \mathbf{y}|||^2}{4}$, where $\mathbf{x}, \mathbf{y} \in \mathbb{W}$.

Definition~\ref{D1} implies that multiplication $\otimes$ of  the basis elements $\mathbb{E}_4$ forms a skew circulated Hankel matrix.
	$$ \begin{array}{c|cccc}
	\otimes & {\mathbf{1}
		} & \mathbf{i}& \mathbf{j} & \mathbf{k}  \\  \hline
	{\mathbf{1}
	}       & {\mathbf{1}
	} & \mathbf{i}& \mathbf{j} & \mathbf{k}  \\
	\mathbf{i}      & \mathbf{i}& \mathbf{j} & \mathbf{k} & -{\mathbf{1}
	}  \\
	\mathbf{j}       & \mathbf{j} & \mathbf{k} & -{\mathbf{1}
	} & -\mathbf{i} \\ 
	\mathbf{k}       & \mathbf{k} & -{\mathbf{1}
	} & -\mathbf{i}& -\mathbf{j}  
	\end{array}  $$
	\end{remark}

Algebra $\mathbb{W}$ belongs to a class of algebras which are generalization of complex numbers in $\mathbb{R}^n$ spaces. For example N. Fleury et al. in \cite{Fleury}, considered an $n$-dimensional commutative algebra generated by n vectors $1,e,,\ldots, e^{n-1}$ where the fundamental element satisfies the basis relation $e^n=-1$. Detailed investigations  for this type of 3-dimensional algebra  were carried out, for example in
\cite{AV,Haluska-Jastrab, LT}.
Algebra  
 $\mathbb{W}$ is a system of  real  numbers generated by element $\mathbf{e}$, which satisfies  the basic relation $\mathbf{e}^4=-1$. 

\begin{theorem}\rm \label{1MJ}  The algebra $\mathbb{W}$   is isomorphic to the algebra $\mathbb{T}$ over $\mathbb{R}$ with the basis $\{1_\mathbb{T}, \mathbf{e}, \mathbf{e}^2, \mathbf{e}^3\}$, where $\mathbf{e}^4=-1_\mathbb{T}$.

$$ \begin{array}{c|cccc}
	\cdot & {\mathbf{1}_\mathbb{T}
	} & \mathbf{e}& \mathbf{e^2} & \mathbf{e^3}  \\  \hline
	{\mathbf{1}_\mathbb{T}}       & {\mathbf{1}_\mathbb{T}} & \mathbf{e}& \mathbf{e^2} & \mathbf{e^3}  \\
	\mathbf{e}      & \mathbf{e}& \mathbf{e^2} & \mathbf{e^3} & -{\mathbf{1}_\mathbb{T}}  \\
	\mathbf{e^2}       & \mathbf{e^2} & \mathbf{e^3} & -{\mathbf{1}_\mathbb{T}} & -\mathbf{e} \\ 
	\mathbf{e^3}       & \mathbf{e^3} & -{\mathbf{1}_\mathbb{T}} & -\mathbf{e}& -\mathbf{e^2}  
	\end{array}  $$
	\end{theorem}
	\proof
The isomorphism on the elements is as follows: $\mathbf{1}\leftrightarrow 1_\mathbb{T}, i \leftrightarrow \mathbf{e}, j\leftrightarrow \mathbf{e}^2, k\leftrightarrow~\mathbf{e}^3$.	\qed

Let $\mathbb{S}$ be an algebra of skew circulant  $4\times 4$ Toeplitz matrices over $\mathbb{R}$, i.e.,

$$\mathbb{S}=\left\{ \left[
\begin{array}{cccc} 
a & b& c &d \\
-d & a& b&c \\ 
-c & -d & a &b \\ -b & -c& -d&a  \end{array} \right] | a, b, c, d \in \mathbb{R}\right\}$$
with  standard matrix operations.
We can establish matrix representation of algebra $\mathbb{W}$ as follows

\begin{theorem}\label{rep1}
Let $\mathbf{x}=(X_\mathbf{1}, X_\mathbf{i}, X_\mathbf{j}, X_\mathbf{k}) \in \mathbb{W}$, where $X_\mathbf{1}, X_\mathbf{i}, X_\mathbf{j}, X_\mathbf{k} \in \mathbb{R}$. Let   $\psi$  be a following  bijection between $\mathbb{W}$ and $\mathbb{S}$:

$$\mathbb{W} \ni  (X_\mathbf{1} X_\mathbf{i}, X_\mathbf{j}, X_\mathbf{k}) \stackrel{\psi}{\mapsto} 
\left[
\begin{array}{cccc} 
X_\mathbf{1}& X_\mathbf{i} & X_\mathbf{j} &X_\mathbf{k} \\
-X_\mathbf{k} & X_\mathbf{1}& X_\mathbf{i} &X_\mathbf{j} \\ -X_\mathbf{j} & -X_\mathbf{k} & X_\mathbf{1}&X_\mathbf{i} \\ -X_\mathbf{i} & -X_\mathbf{j} & -X_\mathbf{k} &X_\mathbf{1}  \end{array} \right] \in \mathbb{S}. $$ 

Then $ {\psi} $ 
 is an isomorphism of algebras $W$ and $\mathbb{S}$.
  \end{theorem}
 
\proof The proof of the   theorem is trivial and left to the reader.  \qed
  
The algebra $ \mathbb{W}$ is isomorphic with some polynomial algebra. The following statement establishes this correspondence.

\begin{theorem} 
Let $\mathbf{x}=(X_\mathbf{1}, X_\mathbf{i}, X_\mathbf{j}, X_\mathbf{k}) \in \mathbb{W}$ where $X_\mathbf{1}, X_\mathbf{i}, X_\mathbf{j}, X_\mathbf{k} \in \mathbb{R}$. Let  $\mathbb{R}[y]$ be an algebra of polynomials  in  one variable $y$ and coefficients from $\mathbb{R},$ and let  $R= \mathbb{R}[y]/(y^4+1)$ be its  homomorphic image.
Then $ \varphi(\mathbf{x})= X_\mathbf{1}+ X_\mathbf{i}\bar{y}+ X_\mathbf{j}\bar{y}^2+ X_\mathbf{k}\bar{y}^3   $ is an isomorphism of  the algebra $\mathbb{W}$ and  the algebra~$R$.
\end{theorem}

\proof
Since every element of the ring $R$ can be expressed uniquely as $a+b\bar{y}+c\bar{y}^2+d\bar{y}^3$ where $a,b,c,d\in\mathbb{R}$ the proof follows from Theorem \ref{1MJ} and simple calculations. 
\qed

\section{$\blacktriangle$-conjugation   in $\mathbb{W}$}\label{Ch3}
\begin{definition}\rm
	Let 
	$\mathbf{x}=  X_\mathbf{1}\ \mathbf{1} \oplus  X_\mathbf{i}\ \mathbf{i}\oplus  X_\mathbf{j} \ \mathbf{j} \oplus  X_\mathbf{k} \ \mathbf{k} \in \mathbb{W}$,
		where $X_\mathbf{1}$, $X_\mathbf{i}$, $X_\mathbf{j}$, $X_\mathbf{k} \in~\mathbb{R}$.
	Then  
	$$(\mathbf{x})^\blacktriangle  \stackrel{def}{=} X_\mathbf{1}\ \mathbf{1} \ominus   X_\mathbf{k} \ \mathbf{i}\ominus X_\mathbf{j} \ \mathbf{j} \ominus  X_\mathbf{i}\ \mathbf{k} \in \mathbb{W}.$$
		 We will call this unary operation to be a \textit{$\blacktriangle$-conjugation.} 
		 	  We will write  simply: $(\mathbf{x})^\blacktriangle = \mathbf{x}^\blacktriangle.$ 	
\end{definition}	
	
Another unary  function is 
$ \mathbf{x} \mapsto \mathbf{x} \otimes \mathbf{x}^\blacktriangle \in \mathbb{W}. $ 
Clearly, 
	 $\mathbf{x}^{\blacktriangle \blacktriangle}= \mathbf{x} $.	

\begin{theorem}\label{conMJ} \rm 
	If	$\mathbf{x} =  X_\mathbf{1}\ \mathbf{1} \oplus X_\mathbf{i}\ \mathbf{i}\oplus X_\mathbf{j} \ \mathbf{j} \oplus X_\mathbf{k} \ \mathbf{k} \in \mathbb{W}$,
then 
$$ \mathbf{x} \otimes \mathbf{x}^\blacktriangle {=} A(\mathbf{x})  \mathbf{1} \oplus B(\mathbf{x})  \theta,$$ 
where
$$A(\mathbf{x}) \stackrel{def}{=}  X_\mathbf{1}^2 + X_\mathbf{i}^2 + X_\mathbf{j}^2 + X_\mathbf{k}^2 ,   $$ 
$$B(\mathbf{x}) \stackrel{def}{=} X_\mathbf{1}X_\mathbf{i}+ X_\mathbf{i}X_\mathbf{j}+  X_\mathbf{j} X_\mathbf{k}  - X_\mathbf{k} X_\mathbf{1},  $$ 
$$\theta \stackrel{def}{=} \mathbf{i}\ominus \mathbf{k} = (0,1, 0, -1),  $$
and $ X_\mathbf{1}, X_\mathbf{i}, X_\mathbf{j}, X_\mathbf{k} \in   \mathbb{R}$. 
\end{theorem}	
\proof
\begin{multline*} \mathbf{x} \otimes \mathbf{x}^\blacktriangle = 
(X_\mathbf{1}^2 + X_\mathbf{i}^2 + X_\mathbf{j}^2 + X_\mathbf{k}^2) \mathbf{1} \\
\oplus (X_\mathbf{1} X_\mathbf{i}+ X_\mathbf{i} X_\mathbf{j}+ X_
\mathbf{j} X_\mathbf{k} ) \mathbf{i}\oplus (X_\mathbf{k} X_\mathbf{1})(-\mathbf{i}) \\
\oplus( X_\mathbf{1}X_\mathbf{j} + X_\mathbf{i}X_\mathbf{k} ) \mathbf{j} \oplus(X_\mathbf{1}X_\mathbf{j} +  X_\mathbf{i}X_\mathbf{k}) (-\mathbf{j}) \\
 \oplus (X_\mathbf{k} X_\mathbf{1})\mathbf{k}\oplus (X_\mathbf{1}X_\mathbf{i} + X_\mathbf{i}X_\mathbf{j} +X_\mathbf{j}X_\mathbf{k})(-\mathbf{k}) \\
=  (X_\mathbf{1}^2 + X_\mathbf{i}^2 + X_\mathbf{j}^2 + X_\mathbf{k}^2 )\mathbf{1}
\oplus  (X_\mathbf{1} X_\mathbf{i}+ X_\mathbf{i} X_\mathbf{j}+  X_\mathbf{j} X_\mathbf{k}  - X_\mathbf{k} X_\mathbf{1})\theta.
\end{multline*}
 
 So,
 $$ \mathbf{x} \otimes \mathbf{x}^\blacktriangle {=} A (\mathbf{x}) \mathbf{1} \oplus B (\mathbf{x}) \theta.
 \QED$$
  
\begin{remark} \rm
From  Theorem~\ref{conMJ} it follows that the function $\mathbf{x} \mapsto \mathbf{x} \otimes \mathbf{x}^\blacktriangle $ 	
maps  elements of the four dimensional space $\mathbb{W}$ to a two-dimensional Euclidean space with a basis 
$\{ \mathbf{1}, \theta\}$.
	\end{remark}

The following lemma shows that  $\blacktriangle$-conjugation is a homomorphism on~$\mathbb{W}$.

\begin{lemma} If $\mathbf{x}= (X_\mathbf{1}, X_\mathbf{i},X_\mathbf{j},  X_\mathbf{k})\in \mathbb{W}$, $ \mathbf{y}= (Y_\mathbf{1},Y_\mathbf{i},Y_\mathbf{j}, Y_\mathbf{k})\in \mathbb{W}$ and $\lambda\in \mathbb{R}$, then 
	$$ (\mathbf{x} \oplus \mathbf{y})^\blacktriangle= \mathbf{x}^\blacktriangle\oplus \mathbf{y}^\blacktriangle, \hskip 1cm 
	 (\mathbf{x} \otimes \mathbf{y})^\blacktriangle= \mathbf{x}^\blacktriangle\otimes \mathbf{y}^\blacktriangle.$$
	and
	$$ (\lambda \mathbf{x})^\blacktriangle= \lambda \mathbf{x}^\blacktriangle$$	
	 \end{lemma}
\proof 
The proof follows from a direct calculations. We omit the details.~\qed

A linear subspace in $\mathbb{W}$, a plane spanned by two vectors $\mathbf{1}$ and $\theta$, has  particular properties from the reduction of operation $\otimes$ to the plane.

\begin{theorem}
Let $\pi^\blacktriangle\stackrel{def}{=} \{ \mathbf{y}= \mathbf{x} \otimes \mathbf{x}^\blacktriangle \mid \mathbf{x} \in  \mathbb{W}\}$  be a plane defined  with  the triple of points  $(\Lambda, \mathbf{1}, \Theta )$, where 
$$\Theta \stackrel{def}{=} \frac{\theta}{\sqrt{2}}.$$
 is a squeezed vector $\theta$.  
Let the plane $\pi^\blacktriangle$ be equipped with the operation $\boxtimes$ which is a reduction of the operation $\otimes$ from the whole space $\mathbb{W}$ to its subplane $\pi^\blacktriangle$. 

Then the operation $\boxtimes$ is defined   by a table 
$$ 
\begin{array}{c|cc}
\boxtimes & \mathbf{1}        & \Theta    \\  \hline
\mathbf{1}       & \mathbf{1}       & \Theta   \\
\Theta & \Theta & \mathbf{1}         
\end{array}
$$
on the plane $\pi^\blacktriangle$ and  the plane $\pi^\blacktriangle$ is  isomorphic to the  hyperbolic (Clifford)  complex  plane with the ,,real unit" $\mathbf{1}$ and the ,,imaginary  unit" $\Theta$, respectively.   
\end{theorem}
\proof For a theory of generalized complex numbers in the plane, see~\cite{Harkin}.

The  points $( \Lambda, \mathbf{1}, \Theta)$   are  non-collinear and mutually different. Therefore, they define a non-degenerate plane $\pi^\blacktriangle=span_\mathbb{R}\{\mathbf{1}_\mathbb{W}, \Theta \}.$
 
To prove that the plane $\pi^\blacktriangle$ equipped with the operation $\boxtimes$  is isomorphic to the hyperbolic complex plane  it is enough  to show that $\Theta \boxtimes \Theta = \mathbf{1}$. Indeed,
 $$\Theta \boxtimes\Theta = \left(\frac{\mathbf{i}\ominus \mathbf{k}}{\sqrt{2}}\right)\boxtimes \left(\frac{\mathbf{i}\ominus \mathbf{k}}{\sqrt{2}}\right) =  \left(\frac{\mathbf{i}\ominus \mathbf{k}}{\sqrt{2}}\right)\otimes \left(\frac{\mathbf{i}\ominus \mathbf{k}}{\sqrt{2}}\right)$$ 
 $$ =\frac{1}{2}(\mathbf{j} \oplus \mathbf{1} \oplus \mathbf{1} \ominus \mathbf{j}) = \frac{\mathbf{j}\ominus \mathbf{j}}{2}  \oplus \mathbf{1} = \Lambda \oplus \mathbf{1} = \mathbf{1}.$$

The operation of addition is satisfied trivially in the plane $\pi^\blacktriangle$ since it is a linear subspace of the Euclidean space  $\mathbb{E}_4$.\qed  	
\section{Values $\mathcal{A}(\mathbf{x}) + \mathcal{B(\mathbf{x})}$ and $\mathcal{A}(\mathbf{x}) - \mathcal{B}(\mathbf{x})$}\label{Ch4}
	
\begin{definition}\label{dpdmMJ}\rm 
	Let $\mathbf{x}=  X_\mathbf{1}\ \mathbf{1}\oplus  X_\mathbf{i}\ \mathbf{i}\oplus  X_\mathbf{j} \ \mathbf{j} \oplus   X_\mathbf{k} \ \mathbf{k} \in \mathbb{W}$, where $X_\mathbf{1}$,  $X_\mathbf{i}$, $X_\mathbf{j}$, $X_\mathbf{k} \in \mathbb{R}$.   
$$\mathcal{A}(\mathbf{x}) \stackrel{def}{=} A(\mathbf{x}) = X_\mathbf{1}^2 + X_\mathbf{i}^2 + X_\mathbf{j}^2 + X_\mathbf{k}^2  ,$$ 
$$\mathcal{B}(\mathbf{x}) {=}  \sqrt {2}{B}(\mathbf{x})  \stackrel{def}{=} \sqrt {2}(X_\mathbf{1} X_\mathbf{i}+ X_\mathbf{i}X_\mathbf{j}+  X_\mathbf{j} X_\mathbf{k}  - X_\mathbf{k} X_\mathbf{1})$$
	and
$$\mathbb{D}_\oplus \stackrel{def}{=} 
\{\mathbf{x} \in  \mathbb{W}| \mathcal{A}(\mathbf{x}) - \mathcal{B}(\mathbf{x})=0\},$$
		$$ \mathbb{D}_\ominus \stackrel{def}{=} \{\mathbf{x} \in  \mathbb{W}|\mathcal{A}(\mathbf{x}) + \mathcal{B}(\mathbf{x}) = 0\} .$$
		\end{definition}	
It can be trivially seen that
\begin{lemma}\label{L2}\ 

If $\mathbf{x} \in \mathbb{D}_\oplus$ and $\mathbf{x} \neq \Lambda$, then $\mathbf{x} \notin \mathbb{D}_\ominus$. 
If $\mathbf{x} \in \mathbb{D}_\ominus$ and $\mathbf{x} \neq \Lambda$, then $\mathbf{x} \notin \mathbb{D}_\oplus$.
\end{lemma}

\begin{lemma}\ The following  equalities are satisfied  for every $\mathbf{x} \in \mathbb{W}:$
 $$4[\mathcal{A}(\mathbf{x}) - \mathcal{B}(\mathbf{x})]
=(1-\sqrt{2})[(X_\mathbf{1}+ X_\mathbf{i})^2+ ( X_\mathbf{i} + X_\mathbf{j})^2+ (X_\mathbf{j} + X_\mathbf{k})^2 + ( X_\mathbf{k} - X_\mathbf{1})^2]$$
$$ + (1+\sqrt{2})[(X_\mathbf{1}- X_\mathbf{i})^2+ ( X_\mathbf{i} - X_\mathbf{j})^2+ (X_\mathbf{j} - X_\mathbf{k})^2 + ( X_\mathbf{k} + X_\mathbf{1})^2];$$
and
 $$4[\mathcal{A}(\mathbf{x}) + \mathcal{B}(\mathbf{x})]
=(1+\sqrt{2})[(X_\mathbf{1}+ X_\mathbf{i})^2+ (X_\mathbf{i} + X_\mathbf{j})^2+ (X_\mathbf{j} + X_\mathbf{k})^2 + 
(X_\mathbf{k} - X_\mathbf{1})^2]$$
$$ + (1-\sqrt{2})[(X_\mathbf{1}- X_\mathbf{i})^2+ ( X_\mathbf{i} - X_\mathbf{j})^2+ (X_\mathbf{j} - X_\mathbf{k})^2 + ( X_\mathbf{k} + X_\mathbf{1})^2].$$
\end{lemma}
\proof
In both  equalities  it is sufficient to simplify the right hand side of equality to obtain the left one. The details are left to the reader.
\qed

\begin{lemma} \label{nomMJ}   { Let $\mathbf{x} =  X_\mathbf{1} \ \mathbf{1} \oplus X_\mathbf{i}\ \mathbf{i}\oplus X_\mathbf{j} \ \mathbf{j} \oplus X_\mathbf{k} \ \mathbf{k} \in \mathbb{W}.$ Then 
$$ \mathcal{A}(\mathbf{x}) - \mathcal{B}(\mathbf{x})= \left(\frac{X_\mathbf{1}}{\sqrt{2}} - X_{\mathbf{i}} + \frac{X_\mathbf{j}}{\sqrt{2}}\right)^2 + \left(\frac{X_\mathbf{1}}{\sqrt{2}} - \frac{X_\mathbf{j}}{\sqrt{2}} + X_{\mathbf{k}}\right)^2 . $$ }
	\end{lemma} 
\proof
 \begin{multline*}\left(\frac{X_\mathbf{1}}{\sqrt{2}} - X_{\mathbf{i}} + \frac{X_\mathbf{j}}{\sqrt{2}}\right)^2 + \left(\frac{X_\mathbf{1}}{\sqrt{2}} - \frac{X_\mathbf{j}}{\sqrt{2}} + X_{\mathbf{k}}\right)^2  \\= 
 \frac{X^2_\mathbf{1}}{2} + X^2_\mathbf{i} +  \frac{X^2_\mathbf{j}}{2} - \sqrt{2} X_\mathbf{1} X_\mathbf{i} +X_\mathbf{1} X_\mathbf{j}  - \sqrt{2} X_\mathbf{j} X_\mathbf{i}   \\ + 
\frac{X^2_\mathbf{1}}{2} + X^2_\mathbf{k} +  \frac{X^2_\mathbf{j}}{2} -  X_\mathbf{1} X_\mathbf{j} +  \sqrt{2} X_\mathbf{1} X_\mathbf{k}  -  \sqrt{2} X_\mathbf{j} X_\mathbf{k}=    \mathcal{A}(\mathbf{x}) - \mathcal{B}(\mathbf{x}). 
 \QED \end{multline*}    
 
\begin{lemma} \label{nopMJ}      { Let $\mathbf{x} =  X_\mathbf{1} \ \mathbf{1} \oplus X_\mathbf{i}\ \mathbf{i}\oplus X_\mathbf{j} \ \mathbf{j} \oplus X_\mathbf{k} \ \mathbf{k} \in \mathbb{W}.$ Then 
$$ \mathcal{A}(\mathbf{x}) +\mathcal{B}(\mathbf{x})=
\left( \frac{X_\mathbf{1}}{\sqrt{2}}  + X_\mathbf{i} + \frac{X_\mathbf{j}}{\sqrt{2}}\right)^2 + \left(\frac{X_\mathbf{1}}{\sqrt{2}} - \frac{X_\mathbf{j}}{\sqrt{2}} - X_{\mathbf{k}}\right)^2 . $$}
 \end{lemma}

 \proof  \begin{multline*}\left(\frac{X_\mathbf{1}}{\sqrt{2}}  + X_\mathbf{i} + \frac{X_\mathbf{j}}{\sqrt{2}}\right)^2 + \left(\frac{X_\mathbf{1}}{\sqrt{2}} - \frac{X_\mathbf{j}}{\sqrt{2}} - X_\mathbf{k}\right)^2 \\  =
   \frac{X^2_\mathbf{1}}{2} + X^2_\mathbf{i} +  \frac{X^2_\mathbf{j}}{2} + \sqrt{2} X_\mathbf{1} X_\mathbf{i} +X_\mathbf{1} X_\mathbf{j}  + \sqrt{2} X_\mathbf{j} X_\mathbf{i} \\ +   
  \frac{X^2_\mathbf{1}}{2}  +\frac{X^2_\mathbf{j}}{2}  
  + X^2_\mathbf{k}  -  X_\mathbf{1}X_\mathbf{j} -  \sqrt{2} X_\mathbf{1}X_\mathbf{k}  +  \sqrt{2} X_\mathbf{j} X_\mathbf{k}= \mathcal{A}(\mathbf{x}) +\mathcal{B}(\mathbf{x})     \QED \end{multline*}
 
 From these two lemmas we obtain explicit  general expressions for $ \mathbb{D}_\oplus$ and~$\mathbb{D}_\ominus$.
    {
 \begin{theorem}\label{kk} Let
  $$ \mathbb{D}_\oplus = \left\{ \left(\alpha, \frac{\alpha + \beta}{\sqrt {2}}, \beta, \frac{-\alpha + \beta}{\sqrt {2}} \right)  \in \mathbb{W}~|~ \alpha, \beta \in \mathbb{R}    \right\}.$$
 Then $$ \mathbb{D}_\oplus =span_\mathbb{R} \{\mathbf{1}_{\mathbb{D}_\oplus}, \mathbf{i}_{\mathbb{D}_\oplus}  \}, $$
 where $$\mathbf{1}_{\mathbb{D}_\oplus} \stackrel{def}{=} \left(\frac{1}{2}, \frac{1}{\sqrt{8}}, 0,  - \frac{1}{\sqrt{8}} \right),  ~ 
  \mathbf{i}_{\mathbb{D}_\oplus} \stackrel{def}{=} \left(0,   \frac{1}{\sqrt{8}}, \frac{1}{2},  \frac{1}{\sqrt{8}} \right),$$
  and 
      $$\mathcal{A}(\mathbf{x})= \mathcal{B}(\mathbf{x}) =2(\alpha^2 + \beta^2) \textrm{ for } \mathbf{x} \in \mathbb{D}_{\oplus}.$$
  Similarly,   let  
  $$ \mathbb{D}_\ominus = \left\{   \left(\gamma, -\frac{\gamma + \delta}{\sqrt {2}}, \delta, \frac{\gamma - \delta}{\sqrt {2}}\right) \in \mathbb{W} ~|~ \gamma, \delta \in \mathbb{R}    \right\} .$$
  Then  $\mathbb{D}_\ominus = span_\mathbb{R} \{ \mathbf{1}_{\mathbb{D}_\ominus}, \mathbf{i}_{\mathbb{D}_\ominus}  \}, $
  where
 $$ \mathbf{1}_{\mathbb{D}_\ominus} \stackrel{def}{=} \left(\frac 12, -\frac{1}{\sqrt{8}}, 0,  \frac{1}{\sqrt{8}} \right), 
 \mathbf{i}_{\mathbb{D}_\ominus} \stackrel{def}{=} \left(0,  - \frac{1}{\sqrt{8}}, \frac{1}{2}, - \frac{1}{\sqrt{8}} \right).$$  
 and 
   $$ \mathcal{A}(\mathbf{x})= - \mathcal{B}(\mathbf{x}) = 2(\gamma^2 +  \delta^2) \textrm{ for } \mathbf{x} \in \mathbb{D}_{\ominus}.$$
 \end{theorem}}
 \proof 
 Let $\mathbf{x}= (X_\mathbf{1}, X_\mathbf{i}, X_\mathbf{j}, X_\mathbf{k})  
 \in  \mathbb{D}_\oplus$. 
 By    { definition of $\mathbb{D}_\oplus$ and } Lemma \ref{nomMJ},
  $$ \frac{X_\mathbf{1}}{\sqrt{2}} - X_{\mathbf{i}} + \frac{X_\mathbf{j}}{\sqrt{2}} =0, \frac{X_\mathbf{1}}{\sqrt{2}} - \frac{X_\mathbf{j}}{\sqrt{2}} + X_{\mathbf{k}}=0.$$
 Putting $\alpha= X_\mathbf{1}$ and $\beta=X_\mathbf{j} $ we obtain
  $$\mathbf{x}
 = \left(\alpha, \frac{\alpha + \beta}{\sqrt {2}}, \beta, \frac{-\alpha + \beta}{\sqrt {2}} \right).
  $$
  
In particular,    substituting $\alpha=\frac 12, \beta=0$ and then   $\alpha=0, \beta =\frac 12$ we obtain:\\
   $$ \mathbb{D}_\oplus =  span_\mathbb{R} \left\{ \left( \frac 12, \frac{1}{\sqrt{8}}, 0,  - \frac{1}{\sqrt{8}} \right),   \left(0,   \frac{1}{\sqrt{8}}, \frac 12, \frac{1}{\sqrt{8}} \right)  \right\}.$$  
    
 The formula 
 $$\mathcal{A}(\mathbf{x})= \mathcal{B}(\mathbf{x}) =2(\alpha^2 + \beta^2), \ \ \alpha, \beta \in \mathbb{R} $$ 
    $\textrm{ for } \mathbf{x} \in \mathbb{D}_{\oplus}$ follows directly from  Definition \ref{dpdmMJ} and a trivial computation.
  
  Analogous description for  $\mathbb{D}_\ominus$ can be obtained by Lemma~\ref{nopMJ}.
   \qed 
 
 \begin{theorem}\label{DSMJ}\ 
 
\begin{itemize}
\item[{(i)}]
 The linear space $\mathbb{D}_\oplus$ with  multiplication~$\otimes$ and with  an identity element $\mathbf{1}_\mathbb{D_\oplus}$ forms a  commutative algebra over field $\mathbb{R}.$ \item[{(ii)}]
 The linear space $\mathbb{D}_\ominus$ with  multiplication~$\otimes$ and with an identity element $\mathbf{1}_\mathbb{D_\ominus}$ forms a  commutative algebra over field $\mathbb{R}.$  
\end{itemize}
\end{theorem}
 \proof
\begin{itemize}
\item[{(i)}] 
 Elements $\mathbf{1}_{\mathbb{D}_\oplus}$ and $\mathbf{i}_{\mathbb{D}_\oplus}$ are basis elements of $\mathbb{D}_\oplus$  and their
  multiplication  table  is as follows:  
$$ 
\begin{array}{c|cc}
\otimes & \mathbf{1}_{\mathbb{D}_\oplus}       & \mathbf{i}_{\mathbb{D}_\oplus}   \\  \hline
\mathbf{1}_{\mathbb{D}_\oplus}      & \mathbf{1}_{\mathbb{D}_\oplus}     & \mathbf{i}_{\mathbb{D}_\oplus}\\
\mathbf{i}_{\mathbb{D}_\oplus} &\mathbf{i}_{\mathbb{D}_\oplus} & -\mathbf{1}_{\mathbb{D}_\oplus}      
\end{array}
$$
Since multiplication in $\mathbb{D}_\oplus$  is determined uniquely by multiplication of basis elements, we obtain that  $\mathbb{D}_\oplus$ is multiplicatively closed set. 
It is evident that  $\mathbf{1}_{\mathbb{D}_\oplus}$ is an identity element in $\mathbb{D}_\oplus$.
\item[{(ii)}]  
The multiplication of basis elements $\mathbf{1}_{\mathbb{D}_\ominus}$ and $\mathbf{i}_{\mathbb{D}_\ominus}$ in $\mathbb{D}_\ominus$ is given by the table

$$ 
\begin{array}{c|cc}
\otimes & \mathbf{1}_{\mathbb{D}_\ominus}       & \mathbf{i}_{\mathbb{D}_\ominus}   \\  \hline
\mathbf{1}_{\mathbb{D}_\ominus}      & \mathbf{1}_{\mathbb{D}_\ominus}     & \mathbf{i}_{\mathbb{D}_\ominus}\\
\mathbf{i}_{\mathbb{D}_\ominus} &\mathbf{i}_{\mathbb{D}_\ominus} & -\mathbf{1}_{\mathbb{D}_\ominus}    
\end{array}
$$

Hence  $\mathbb{D}_\ominus$ is multiplicatively closed set. It is clear that $\mathbf{1}_{\mathbb{D}_\ominus}$ is an identity element in~$\mathbb{D}_\ominus.$  \qed \end{itemize}

\begin{lemma} \label{L6}\rm
Let $\mathbf{x}\in \mathbb{D}_\oplus$ and $\mathbf{y} \in  \mathbb{D}_{\ominus}$. Then $\mathbf{x}\otimes\mathbf{y} =\Lambda$. 
\end{lemma}
\proof 
There exist $\alpha, \beta \in \mathbb{R}$ and $    {\gamma, \delta} \in \mathbb{R}$ such that 
$$ \mathbf{x} = 
\left( \alpha, \frac{\alpha+ \beta}{\sqrt{2}}, \beta, \frac{-\alpha+ \beta}{\sqrt{2}}\right), \mathbf{y} = 
   {\left(\gamma, - \frac{\gamma + \delta}{\sqrt{2}}, \delta,  \frac{\gamma- \delta}{\sqrt{2}}\right). }$$
By a simply verification by Definition~1, we obtain $\mathbf{x}\otimes\mathbf{y} =\Lambda$. 
\qed

\begin{theorem}\label{TH8}
The spaces $\mathbb{D}_\oplus$, $\mathbb{D}_\ominus$ are ideals in 
$\mathbb{W}$, and  $\mathbb{W}=\mathbb{D}_\oplus \times \mathbb{D}_\ominus,$ as a direct sum of  ideals. Moreover, $\mathbb{D}_\oplus$ is  isomorphic to $\mathbb{C}$ and also  $\mathbb{D}_\ominus$  is isomorphic to $\mathbb{C}$.   Therefore, $\mathbb{W}$ is isomorphic to $\mathbb{C}\times \mathbb{C}$.
\end{theorem}
\proof
Since
        $$\det \left[ 
        \begin{array}{cccc}
        \frac 12 &\frac{1}{\sqrt {8}} & 0 &- \frac{1}{\sqrt {8}}\\
        0&  \frac{1}{\sqrt {8}}& \frac 12& \frac{1}{\sqrt {8}} \\
       \frac 12&-  \frac{1}{\sqrt {8}}&0&\frac{1}{\sqrt {8}} \\
       0& -\frac{1}{\sqrt {8}}& \frac 12& - \frac{1}{\sqrt {8}}\end{array}
       \right] \neq 0, $$ 
this proves that    $\Lambda=\mathbb{D}_\oplus \cap \mathbb{D}_\ominus$ and    $\mathbb{W}=\mathbb{D}_\oplus + \mathbb{D}_\ominus$, as a direct sum of linear spaces.

By  Lemma~\ref{L6} and Theorem \ref{DSMJ},  $\mathbb{D}_\oplus$ and $\mathbb{D}_\ominus$ are ideals in 
$\mathbb{W}$ and $\mathbb{W}=\mathbb{D}_\oplus \times \mathbb{D}_\ominus$. 
It remains to prove that  $\mathbb{D}_\oplus$ is isomorphic to $\mathbb{C}$ and analogously that  $\mathbb{D}_\ominus$ is isomorphic to $\mathbb{C}.$

A trivial verification shows that  $\left( 0, -\frac{1}{\sqrt{8}}, -\frac{1}{2}, - \frac{1}{\sqrt{8}} \right)$ is a (Gaussian) "imaginary unit" in $\mathbb{D}_\oplus$ and $\left( 0, \frac{1}{\sqrt{8}}, -\frac{1}{2},  \frac{1}{\sqrt{8}} \right)$ is a (Gaussian) "imaginary unit" in $\mathbb{D}_\ominus$, respectively.  
Indeed,  under the notation in the proof of Theorem \ref{DSMJ}, we have :
$$ \left( 0, -\frac{1}{\sqrt{8}}, -\frac{1}{2}, - \frac{1}{\sqrt{8}} \right)^2= \left(-\mathbf{i}_{\mathbb{D}_\oplus}\right)\otimes \left(-\mathbf{i}_{\mathbb{D}_\oplus}\right) = 
-\mathbf{1}_{\mathbb{D}_\oplus}$$
and
$$\left( 0, \frac{1}{\sqrt{8}}, -\frac{1}{2}, \frac{1}{\sqrt{8}} \right)^2=   \left(-\mathbf{i}_{\mathbb{D}_\ominus}\right) \otimes \left(-\mathbf{i}_{\mathbb{D}_\ominus}\right)= -\mathbf{1}_\mathbb{D_\ominus}.$$ 
This completes the proof.
\qed 

In other words, the theorem states  that  $\mathbf{x} \in  \mathbb{W}$ is a divisor of zero if and only if $\mathbf{x} \in  \mathbb{D}_\oplus$ or $\mathbf{x} \in  \mathbb{D}_\ominus$. Using isomorphism  $\mathbb{D}_\oplus$  to $\mathbb{C}$ and   $\mathbb{D}_\ominus$  to  $\mathbb{C}$ one can determine the inverse elements for each invertible element in $\mathbb{W}$. But the calculations are quite complicated. 

\section{Alternative method of finding   inverse elements}

In this section  we find an easier method to compute  inverse elements in $\mathbb{W}$ using the skew circulated structure of the definition of the operation $\otimes$.

\begin{theorem}\label{TH10}
Let $\mathbf{x}=  X_\mathbf{1} \mathbf{1} \oplus X_\mathbf{i} \mathbf{i}  \oplus X_\mathbf{j} \mathbf{j} \oplus   X_\mathbf{k}\mathbf{k}$,
	where $ X_\mathbf{1}, X_\mathbf{i}, X_\mathbf{j}$, $X_\mathbf{k} \in \mathbb{R}$. 
	
Let $\mathbf{x} \not\in  \mathbb{D}_\oplus$  and   $\mathbf{x} \not\in  \mathbb{D}_\ominus$.   
Then the multiplicative inverse  to the element  $\mathbf{x}$ exists and it is unique. 
The  element  $\mathbf{x^{-1}}$ is  defined as follows:
	$$ \mathbf{x}^{-\mathbf{1}}  
	\stackrel{def}{=}  \frac{\mathbf{x}^\blacktriangle\otimes [\mathcal{A}(\mathbf{x})\ \mathbf{1} \ominus \mathcal{B}(\mathbf{x}) \ \Theta]}{[\mathcal{A}(\mathbf{x}) + \mathcal{B}(\mathbf{x})][\mathcal{A}(\mathbf{x}) - \mathcal{B}(\mathbf{x})]},$$ 
 where $\mathbf{x}\otimes\mathbf{x}^{-\mathbf{1}}
	= \mathbf{1}$. 
\end{theorem} 
\proof
 Since $\mathbf{x} \not\in  \mathbb{D}_\oplus$  and $\mathbf{x} \not\in  \mathbb{D}_\ominus$ then $[\mathcal{A}(\mathbf{x}) +  \mathcal{B}(\mathbf{x}) ][\mathcal{A}(\mathbf{x}) - \mathcal{B}(\mathbf{x}) ]\neq 0.$

By Theorem \ref{conMJ},
\begin{multline*}\mathbf{x}\otimes\mathbf{x}^{-{\mathbf{1}}} 
	 = \mathbf{x} \otimes \frac{\mathbf{x}^\blacktriangle\otimes[\mathcal{A}(\mathbf{x})\ \mathbf{1} \ominus \mathcal{B}(\mathbf{x}) \ \Theta]}{[\mathcal{A}(\mathbf{x}) +  \mathcal{B}(\mathbf{x}) ] [\mathcal{A}(\mathbf{x}) - \mathcal{B}(\mathbf{x}) ]} \\ 
 =  \frac{ (\mathbf{x} \otimes \mathbf{x}^\blacktriangle)\otimes[\mathcal{A}(\mathbf{x})\ \mathbf{1} \ominus \mathcal{B}(\mathbf{x}) \ \Theta]}
	 {[\mathcal{A}(\mathbf{x}) +  \mathcal{B}(\mathbf{x}) ] [\mathcal{A}(\mathbf{x}) - \mathcal{B}(\mathbf{x}) ]} 	 
\\	 = \frac{ [\mathcal{A}(\mathbf{x})\ \mathbf{1} \oplus \mathcal{B}(\mathbf{x}) \ \Theta]
	 \otimes[\mathcal{A}(\mathbf{x})\ \mathbf{1} \ominus \mathcal{B}(\mathbf{x}) \ \Theta]}{[\mathcal{A}(\mathbf{x}) +  \mathcal{B}(\mathbf{x}) ] [\mathcal{A}(\mathbf{x}) - \mathcal{B}(\mathbf{x}) ]} \\
	  =\frac{[ \mathcal{A}^2(\mathbf{x})-\mathcal{B}^2(\mathbf{x})] \  \mathbf{1}}
	  {[\mathcal{A}(\mathbf{x}) +  \mathcal{B}(\mathbf{x}) ] [\mathcal{A}(\mathbf{x}) - \mathcal{B}(\mathbf{x}) ]}=
	   \mathbf{1}. \end{multline*}

                  The uniqueness of $\mathbf{x}^{-\mathbf{1}}$ is obvious. \qed
 
\section{Topology on $\mathbb{W}$}
\begin{definition} Under the previous notation of $\mathcal{A}   {(\mathbf{x})}$ and $\mathcal{B}   {(\mathbf{x})}$, let us denote by
$$ \| \cdot\|_   {\ominus} \stackrel{def}{=} \sqrt{\mathcal{A}( \cdot ) - \mathcal{B}(\cdot)}:    {\mathbb{D}_\ominus }\to [0,\infty)$$ 
and
$$ \| \cdot \|_   {\oplus} \stackrel{def}{=} \sqrt{\mathcal{A}(\cdot ) + \mathcal{B}(\cdot )}:     {\mathbb{D}_\oplus} \to [0,\infty).$$
\end{definition}

\begin{theorem}  The function $\| \cdot \|_   {\ominus}:    {\mathbb{D}_\ominus} \to [0, + \infty)$  is  a norm.
   \end{theorem}
\proof   Let us test the norm  definition.  

Let $\mathbf{x}=(X_\mathbf{1}, X_\mathbf{i}, X_\mathbf{j}, X_\mathbf{k}) \in    {\mathbb{D}_\ominus}$.
  The following properties 
\begin{itemize}
\item[{(i)}]
$\|\mathbf{x}\|_   {\ominus} = 0 \Leftrightarrow \mathbf{x} =\Lambda$;  
\item[{(ii)}]
 $\|\mathbf{x} \|_   {\ominus} \geq 0$;
\item[{(iii)}]   $\|\alpha \mathbf{x} \|_   {\ominus}  =|\alpha| \|\mathbf{x} \|_   {\ominus} $ for every $\alpha \in \mathbb{R}$
\end{itemize}
hold true    {  for every $\mathbf{x}\in \mathbb{D}_\ominus$ by Theorem \ref{kk}   and simple calculations.}  

Let us prove the triangle inequality. To do this, let another element $ \mathbf{y} =(Y_\mathbf{1}, Y_\mathbf{i}, Y_\mathbf{j}, Y_\mathbf{k})\in    {\mathbb{D}_\ominus}$.    {Further,}  
let us  denote by $$\|\mathbf{x}\|_   {\ominus}  = \sqrt{r_\mathbf{x}^2 + s_\mathbf{x}^2}, 
\hskip 1cm r_\mathbf{x} =  \frac{X_\mathbf{1}}{\sqrt{2}} - X_\mathbf{i} + \frac{X_\mathbf{j}}{2}, 
\hskip 1cm 
s_\mathbf{x} =  \frac{X_\mathbf{1}}{\sqrt{2}} - \frac{X_\mathbf{j}}{\sqrt{2}} + X_\mathbf{k},$$ 
$$\|\mathbf{y} \|_   {\ominus} =\sqrt{r_\mathbf{x}^2 + s_\mathbf{x}^2}, 
\hskip 1cm r_\mathbf{y} =  \frac{Y_\mathbf{1}}{\sqrt{2}} - Y_\mathbf{i} + \frac{Y_\mathbf{j}}{\sqrt{2}}, 
\hskip 1cm s_\mathbf{y}=  \frac{Y_\mathbf{1}}{\sqrt{2}} - \frac{Y_\mathbf{j}}{\sqrt{2}} + Y_\mathbf{k},$$
 and
$$\|\mathbf{x \oplus y}\|_   {\ominus}  = \sqrt{r_\mathbf{x \oplus y}^2 + s_\mathbf{x \oplus y}^2},$$ 
where  
$$r_\mathbf{x \oplus y} =  \frac{X_\mathbf{1}+ Y_\mathbf{1}}{\sqrt{2}} - (X_\mathbf{i} + Y_\mathbf{i}) + \frac{X_\mathbf{j} + Y_\mathbf{j}}{2}, 
\hskip 1cm 
s_\mathbf{x \oplus y}=  \frac{X_\mathbf{1}+ Y_\mathbf{1}}{\sqrt{2}} - \frac{X_\mathbf{j} + Y_\mathbf{j}}{\sqrt{2}} + (X_\mathbf{k}+Y_\mathbf{k}).$$ 

We have:
$$r_\mathbf{x\oplus y} =  \frac{X_\mathbf{1}+ Y_\mathbf{1}}{\sqrt{2}} - (X_\mathbf{i} + Y_\mathbf{i}) + \frac{X_\mathbf{j}+ Y_\mathbf{j}}{2} = r_\mathbf{x} +r_\mathbf{y},$$ 
and
 $$ s_\mathbf{x \oplus y}=  \frac{X_\mathbf{1} + Y_\mathbf{1}}{\sqrt{2}} - \frac{X_\mathbf{j} + Y_\mathbf{j}}{\sqrt{2}} + (X_\mathbf{k}+Y_\mathbf{k}   {)} = s_\mathbf{x} +s_\mathbf{y}.$$
 
 The inequality
 $$ 0 \leq (r_\mathbf{x} s_\mathbf{y} - s_\mathbf{x} r_\mathbf{y})^2$$
  implies  
 $$r^2_\mathbf{x} r^2_\mathbf{y} + s^2_\mathbf{x}s^2_\mathbf{y} + 2r_\mathbf{x} r_\mathbf{y} s_\mathbf{x} s_\mathbf{y} \leq r^2_\mathbf{x} r^2_\mathbf{y} + r^2_\mathbf{x} s^2_\mathbf{y} + s^2_\mathbf{x} r^2_\mathbf{y} + s^2_\mathbf{x} s^2_\mathbf{y} $$
 which implies
 $$ (r_\mathbf{x} r_\mathbf{y} + s_\mathbf{x} s_\mathbf{y})^2 \leq (r_\mathbf{x}^2 + s_\mathbf{x}^2)(r_\mathbf{y}^2 + s_\mathbf{y}^2)$$
 hence we have
 $$ r^2_\mathbf{x} + r^2_\mathbf{y} + 2 r_\mathbf{x} r_\mathbf{y} + s^2_\mathbf{x} + s^2_\mathbf{y} + 2 s_\mathbf{x} s_\mathbf{y} \leq r^2_\mathbf{x} + s^2_\mathbf{x} + r^2_\mathbf{y} + s^2_\mathbf{y} + 2 \sqrt{(r^2_\mathbf{x} + s^2_\mathbf{x} )(r^2_\mathbf{y} + s^2_\mathbf{y})}$$
 and, finally,
 $$\sqrt{(r_\mathbf{x} + r_\mathbf{y})^2 + (s_\mathbf{x} + s_\mathbf{y})^2} \leq  \sqrt{r_\mathbf{x}^2 + s_\mathbf{x}^2} +  \sqrt{r^2_\mathbf{y} + s^2_\mathbf{y}} .$$
In other words,
 $$\|\mathbf{x} \oplus \mathbf{y} \|_   {\ominus} \leq \|\mathbf{x}\|_   {\ominus}  + \| \mathbf{y}\|_   {\ominus}.  \QED $$






Analogously, we can consider  the space    {$\mathbb{D}_\oplus$} and obtain the following theorem.

\begin{theorem} The function $\| \cdot \|   {_\oplus: \mathbb{D}_{\oplus}} \to [0, + \infty)$  is  a norm.
\end{theorem}
\proof
   Let $\mathbf{x}=(X_\mathbf{1}, X_\mathbf{i}, X_\mathbf{j}, X_\mathbf{k}) \in \mathbb{D}_\oplus$.
The following properties 
 \begin{itemize}
 \item[{(i)}]
 $\|\mathbf{x}\|_   {\oplus} = 0 \Leftrightarrow \mathbf{x} =\Lambda$;  
 \item[{(ii)}]
  $\|\mathbf{x} \|_   {\oplus} \geq 0$;
 \item[{(iii)}]   $\|\alpha \mathbf{x} \|_\oplus =|\alpha| \|\mathbf{x} \|_   {\oplus}$ for every $\alpha \in \mathbb{R}$
 \end{itemize}
 
 hold true for every $\mathbf{x}\in \mathbb{D}_\oplus$ by Theorem \ref{kk} and simple calculations.

 Let us prove the triangle inequality. To do this, let another element $ \mathbf{y} =(Y_\mathbf{1}, Y_\mathbf{i}, Y_\mathbf{j}, Y_\mathbf{k})\in \mathbb{D}_   {_\oplus}$.    {Further,}
 
For the triangle inequality, i.e.,  
$$\|\mathbf{x} \oplus \mathbf{y} \|_   {_\oplus} \leq \|\mathbf{x}\|_    {_\oplus}+ \| \mathbf{y}\|   {_\oplus}$$
we put 
$$\|\mathbf{x}\|_   {_\oplus} = \sqrt{(r'_\mathbf{x})^2 + (s'_\mathbf{x})^2}, 
\hskip 1cm r'_\mathbf{x} =  \frac{X_\mathbf{1}}{\sqrt{2}} + X_\mathbf{i} + \frac{X_{\mathbf{j}}}{2}, 
\hskip 1cm s'_\mathbf{x} =  \frac{X_\mathbf{1}}{\sqrt{2}} - \frac{X_\mathbf{j}}{\sqrt{2}} - X_{\mathbf{k}},$$ 
$$\|\mathbf{y}\|_   {_\oplus} = \sqrt{(r'_\mathbf{x})^2 + (s'_\mathbf{x})^2}, \hskip 1cm r'_\mathbf{y} =  \frac{Y_\mathbf{1}}{\sqrt{2}} + Y_\mathbf{i} + \frac{Y_{\mathbf{j}}}{\sqrt{2}}, \hskip 1cm s'_\mathbf{y}=  \frac{Y_\mathbf{1}}{\sqrt{2}} - \frac{Y_\mathbf{j}}{\sqrt{2}} - Y_{\mathbf{k}}.$$ 

The rest of proving    { the triangle inequality} is similar as in the previous proof.  \qed

Since a direct sum of two normed spaces  is a normed space, we define:

\begin{definition} Let $\mathbf{x} \in \mathbb{W} = \mathbb{D}_\oplus \times \mathbb{D}_\ominus$ and  $\mathbf{x}_\oplus \in \mathbb{D}_\oplus$ and $\mathbf{x}_\ominus \in \mathbb{D}_\ominus$ be two elements  such that $\mathbf{x} = \mathbf{x}_\oplus \oplus \mathbf{x}_\ominus$.
The function $$\| \mathbf{x} \| \stackrel{def}{=}  \| \mathbf{x}_\oplus \|_\oplus \oplus \| \mathbf{x}_\ominus \|_\ominus: \mathbb{W} = \mathbb{D}_\oplus \times \mathbb{D}_\ominus  \to [0, +\infty)$$ is a norm equivalent to $||| \cdot |||$.
\end{definition}

\begin{remark}
 This norm is equivalent to the classical Euclidean norm on~$\mathbb{E}_4$ since all norms on a finite $n$ dimensional space are mutually equivalent, $n~\in~\mathbb{N}$.  
 \end{remark}


\end{document}